\documentclass[]{article}

\usepackage{amsmath}
\usepackage{amssymb}
\usepackage{graphicx}

\usepackage{algorithm}
\usepackage{algpseudocode}
\usepackage{adjustbox}
\usepackage{xcolor}
\usepackage{newverbs}
\usepackage{wrapfig}
\usepackage{url}

\usepackage{lscape}

\newlength{\maxwidth}
\newcommand{\algalign}[2]
{\makebox[\maxwidth][r]{$#1{}$}${}#2$}

\newtheorem{theorem}{Theorem}
\newtheorem{lemma}{Lemma}

\newcommand{\surf}{\mathcal{M}_h}
\newcommand{\vfsurf}{\mathcal{X}_h}
\newcommand{\lagrange}{\mathcal{L}}
\newcommand{\crouzeix}{\mathcal{F}}
\newcommand{\lagrangebd}{\mathcal{L}_0}
\newcommand{\crouzeixbd}{\mathcal{F}_0}
\newcommand{\nedelec}{\mathcal{N}}
\newcommand{\nedelecbd}{\mathcal{N}_0}

\newcommand{\textNedelec}{\mbox{N\'{e}d\'{e}lec }}
\newcommand{\smoothsurf}{\mathcal{M}}
\newcommand{\curl}{\mbox{curl}}
\newcommand{\divergence}{\mbox{div}}
\newcommand{\dirichlet}{\mathcal{H}_{h,D}}
\newcommand{\neumann}{\mathcal{H}_{h,N}}
\newcommand{\harmonic}{\mathcal{H}_{h}}

\newcommand{\xnedelecbd}{X_{\mbox{\curl}(\nedelec_0)}}

\newcommand{\xcourzeix}{X_{\nabla\crouzeix}}
\newcommand{\xdirichlet}{X_{\dirichlet}}
\newcommand{\xneumann}{X_{\neumann}}
\newcommand{\xharmcoex}{X_{\harmonic\cap\nabla\crouzeix}}
\newcommand{\xcentre}{X_{\centre}}
\newcommand{\centre}{\mbox{curl}(\nedelec)\cap\nabla\crouzeix}

\newcommand{\fenics}{\texttt{FEniCS }}

\newcommand{\paraview}{\texttt{ParaView }}

\title{A Consistent Discrete 3D Hodge-type Decomposition: implementation and practical evaluation}
\author{Faniry H. Razafindrazaka, Konstantin Poelke, Konrad Polthier,\\ and Leonid Goubergrits}
\date{}
\begin{document}

\maketitle

\begin{abstract}
	The Hodge decomposition provides a very powerful mathematical method for the analysis of 2D and 3D vector fields. 
It states roughly that any vector field can be $L^2$-orthogonally decomposed into a curl-free, divergence-free, and a harmonic field. 
The harmonic field itself can be further decomposed into three components, 
two of which are closely tied to the topology of the underlying domain. 
For practical computations it is desirable to find a discretization which preserves as many aspects inherent to the smooth theory as possible while at the same time remains computationally tractable, in particular on large-sized models. 
The correctness and convergence of such a discretization depends strongly on the choice of ansatz spaces defined on the surface or volumetric mesh to approximate infinite dimensional subspaces. 
This paper presents a consistent discretization  of Hodge-type decomposition for piecewise constant vector fields on volumetric meshes. 
Our approach is based on a careful interplay between edge-based \textNedelec elements and face-based Crouzeix-Raviart elements resulting in a very simple formulation. The method is stable under noisy vector field and mesh resolution, and has a good performance for large sized models. 
We give pseudocodes for a possible implementation of the method together with some insights on how the Hodge decomposition  could answer some central question in computational fluid.
\end{abstract}

\section{Introduction}
\label{sec:intro}
The study of discrete Hodge-type decompositions for 3D tetrahedral meshes poses a challenging problem in vector field analysis. 
Due to its heavy usage in flow analysis of vector fields modeling e.g. blood flow, wind flow or gas propagation it is a central topic of research.  
In addition to the rich theoretical aspects such as the coupling to underlying topological features of the mesh it is important for numerical computations to find a reliable discretization that is simple and efficient to compute, even on large meshes.

Given a vector field on a bounded domain in 3-space, it is a common question to find out whether this field is the gradient of a potential function, the curl of another vector field, or if it is both divergence-free and curl-free  at the same time---a so-called harmonic field.
An answer to this question is provided by the Hodge decomposition theorem which states that the space of vector fields can be decomposed into five mutually $L^2$-orthogonal subspaces. 

More precisely,
let $\smoothsurf$ be a Riemannian 3-manifold in $\mathbb{R}^3$ with smooth boundary $\partial \smoothsurf$. Let $\mathcal{X} = \mathcal{X}(\smoothsurf)$ denote the space of all smooth vector fields on $\smoothsurf$, equipped with the $L^2$-product
\[
\left<X,Y\right>_{L^2} = \int_{\smoothsurf}\left<X,Y\right>_{\mathbb{R}^3}\mbox{d}v,\,\,\mbox{for }X,Y\in \mathcal{X}(\smoothsurf)
\]	
where d$v$ is a volume form on $\mathcal{M}$. 

Following Cantarella's notation~\cite{cantarella_2002}, 
the space $\mathcal{X}$ is the $L^2$-orthogonal sum of five subspaces
\begin{equation}
\label{eq:hodge}
\mathcal{X} = \mbox{FK}\oplus\mbox{GG}\oplus\mbox{CG}\oplus\mbox{HG}\oplus\mbox{HK}
\end{equation}
defined as 
\begin{align*}
\mbox{FK} &= \{\mbox{curl}(Y)\,\,|\,\,\mbox{$Y$ is normal to }\partial\smoothsurf\}\\
\mbox{GG} &= \{\nabla \varphi \,\,|\,\,\varphi_{|\partial \smoothsurf} = 0\}\\
\mbox{CG} &= \{\nabla \varphi \,\,|\,\,\exists Y\in\mathcal{X} \mbox{ such that } \nabla\varphi = \curl(Y)\}\\
\mbox{HG} &= \{X = \nabla \varphi \,\,|\,\,\divergence(X) = 0\mbox{ and $X$ is normal to }\partial\smoothsurf\}\\
\mbox{HK} &= \lbrace X = \curl(Y)\,\,|\,\,\curl(X) = 0\\
&\quad\quad\quad\quad\mbox{ and $X$ is tangential to }\partial\smoothsurf\rbrace.
\end{align*}

\begin{figure}[!h]
	\centering
	\includegraphics[width=\linewidth]{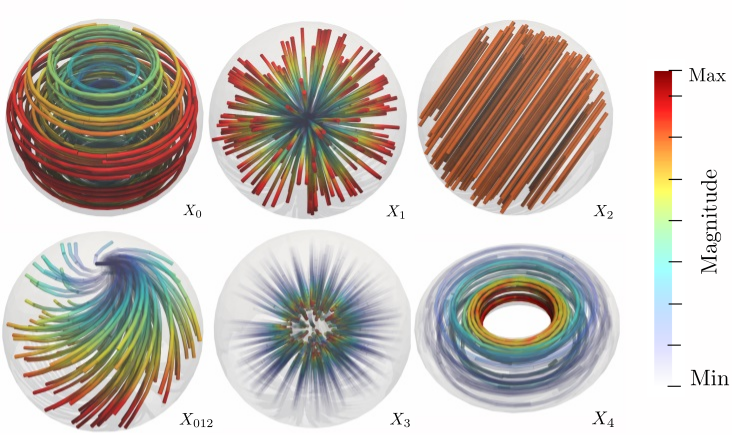}
	\caption{$X_0$-$X_4$ are examples of vector fields belonging to the five subspaces FK, GG, CG, HG, and HK. 
		$X_{012}$ is a superposition of the fields $X_0$-$X_2$ and therefore in the vector space sum FK $+$ GG $+$ CG. 
		Red and opaque encode high magnitude while transparent and blue encode low magnitude.}
	\label{fig:example}
\end{figure}

The elements in FK are called \textit{fluxless knots}.
They are vector fields tangential to $\partial\smoothsurf$ and at 
the same time the curl of a vector field normal to the boundary. 
Elements in GG are called \textit{grounded gradients}. 
These are all vector fields that are gradients of some potential functions $\varphi$ such that $\varphi$ is zero at $\partial \mathcal{M}$. 
The elements in CG are called \textit{curly gradients}.
As the name suggests, these are fields which are simultaneously the gradient of a potential function and the curl of a vector field. 
The spaces HG (\textit{harmonic knots}) and HK (\textit{harmonic gradients}) are the spaces of harmonic fields which are both curl-free and divergence-free, 
and are normal respectively tangential to the boundary $\partial \smoothsurf$. 
In fact, these two spaces HG and HK are always finite-dimensional, and their dimensions depend solely on the topology of $\mathcal{M}$.
The following analytical vector fields  belong to these spaces by choosing the correct domain $\mathcal{M}$.
\begin{align}
X_0 &= (y,-x,0) = \curl\left((-zx,-zy,-z^2)\right)\nonumber\\
X_1 &= (x,y,z) = \frac{1}{2}\nabla(x^2+y^2+z^2-1)\nonumber\\
X_2 &= -\frac{1}{2}(1,1,1)\label{eq:vfsmooth}\\
X_3 &= \frac{1}{(x^2+y^2+z^2)^{3/2}}(x,y,z)\nonumber\\ 
X_4 &= \frac{1}{x^2+y^2}(y,-x,0)\nonumber
\end{align}

By choosing $\mathcal{M}=B(0,1)$ to be the unit ball, it is a simple calculation to check that $X_0\in \mbox{FK}$, $X_1\in \mbox{GG}$ and $X_2\in \mbox{CG}$. 
The vector field $X_3$ belongs to HG for $\mathcal{M} = B(0,1)-B(0, \epsilon)$ for some $0< \epsilon < 1$, a solid unit ball with a small cavity centered at the origin. 
Finally, $X_4\in \mbox{HK}$ lives on a solid torus. 
Figure~\ref{fig:example} shows visualizations of these vector fields on their respective domains. 
Notice that some of these spaces might be trivial, 
depending on the topology of the underlying manifold. 
For example, on a ball-like manifold there are neither harmonic gradients nor harmonic knots.



\subsection{Problem Statement} 
Finding a discrete approximation for the Hodge decomposition (\ref{eq:hodge}) is a delicate issue, as 
it is desirable to reflect the rich structure present in the smooth theory as best as possible, in particular if one aims for a discretization that is as simple as possible.  
For instance, the $L^2$-orthogonality statements should hold for all corresponding discrete spaces, too. 
Second, all topological invariant quantities related to the cohomology of the domain need to be reflected in the corresponding subspaces. 
More precisely, the discrete versions of HG and HK should have the correct dimension in order for the decomposition to be globally consistent. 
Finally, under the assumption that the mesh converges metrically to a smooth manifold, the discrete decomposition should of course converge to the smooth decomposition in some sense.

\subsection{Main Contribution}
We review a set of consistent discrete Hodge-type decompositions on simplicial 3-manifolds embedded in $\mathbb{R}^3$ which fulfill $L^2$-orthogonality, preserve the cohomological dimensions of the harmonic spaces, and converge to the smooth setting in the limit. 
While most of the theoretical proofs are published in~\cite{poelke2017}, we provide in this paper numerical evidence and practical evaluation of the method on simulated vector fields. 
The analysis  proves reliability for large-sized models (millions of tetrahedra), stability under mesh refinement, and robustness against noisy vector fields.  
Our formulation is based  directly on vector proxies and affine spaces.  
Differential forms are only present as side remarks to highlight the  theoretical analogies.	
In addition we study several applications of the decomposition in 
computational fluid dynamics (CFD) of some artificial as well as a patient specific aorta model, giving rise to a new understanding of the blood flow properties under anatomical change as well as variable boundary conditions. 
Our method could be equally applied to 3D velocity fields derived from computational fluid dynamics as well as clinically measured flow fields by 4D flow MRI provided that the domain is tetrahedralized.

\subsection{Related Work}


A comprehensive treatment of Hodge decompositions up to the four-term Hodge-Morrey-Friedrichs decomposition for differential forms is given by Schwarz in~\cite{schwarz_hodge_decomp}. 
A five-term decomposition for vector fields on domains in $\mathbb{R}^3$ is derived in the expository article~\cite{cantarella_2002}. The distinction between fields representing inner and boundary
cohomology and the question of orthogonality between the spaces of Neumann and Dirichlet fields is first published by Shonkwiler in~\cite{shonkwiler_thesis,shonkwiler_2013}.
The focus on picewise constant vector fields (PCVFs) on simplicial surface meshes in modern geometry processing and their theoretical description goes back at least to Polthier and Preuss~\cite{polthier_preuss} who used them as a concept for analysis and decomposition of vector fields, with a convergence proof on closed surfaces given by Wardetzky in~\cite{wardetzky2006thesis}. 
Its generalization to surfaces with boundary is given by Poelke and Polthier in~\cite{Poelke2016126} 
and the theoretical foundations and discretizations discussed in this article are published in~\cite{poelke2017}, 
whereas a first extension to three-dimensional domains based purely on Lagrange ansatz spaces is given in~\cite{Tong_2006}. 

\paragraph*{\textbf{On finite element discretizations}}	
In the 1970s and 1980s, Raviart and Thomas~\cite{raviart1977mixed}
and \textNedelec \cite{nedelec_1980} developed families of finite elements for the computational treatment
of Navier-Stokes- and Maxwell-type problems, and Bossavit \cite{bossavit1988whitney} later put these finite
elements in the context of Whitney forms in two and three dimensions.
The finite element exterior calculus developed by Arnold et al.~\cite{arnold_2006,arnold2010finite} provides
a unifying framework for mixed problems involving more than one ansatz space for finite elements discretization. 
The definition of a discrete curl and divergence for PCVFs as those elements that enforce
Green's formula to hold true is driven by a paradigm inherent to mimetic discretization methods~\cite{bochev_hyman}. Here, the definition of discrete operators and objects is primarily driven by the attempt to mimic properties from the smooth world such as conservation laws, Green's formula or the complex property of operators. 

\paragraph*{\textbf{On other approaches}}
Finally, it should be noted that there are many different definitions and approaches for the
computation of a discrete Hodge decomposition. For instance, ~\cite{bochev_hyman} discuss a decomposition in the spectral domain after Fourier transformation, and~\cite{petronetto2010meshless} propose a meshless decomposition on point set samples in $\mathbb{R}^3$. See also~\cite{bhatia_survey} for a survey on various discrete Hodge decompositions. 
For more application-oriented approaches see e.g. \cite{Bhatia2014}.
The recent work by Zhao et al.~\cite{Zhao2019} formulates a decomposition similar to the one given in \cite{poelke2017}, this time within the framework of discrete exterior calculus. 

\section{Background}
In this section we review some theory and terminology necessary to understand our decomposition based on finite element spaces. 
All our volumetric meshes are assumed to be discretized by tetrahedra. 
As a consequence there is a wide range of classical finite element function spaces available on such meshes, e.g. Lagrange, Crouzeix-Raviart, or \textNedelec spaces.  
The choice of proper function spaces is critical in order to obtain consistent discretizations for Hodge-type decompositions. 
In the following we introduce finite-dimensional spaces on tetrahedral meshes approximating each space FK, GG, CG, HG, and HK introduced in equation~\ref{eq:hodge}.

\subsection{Simplicial Complexes}
\begin{figure}[!ht]
	\centering
	\includegraphics[width=\linewidth]{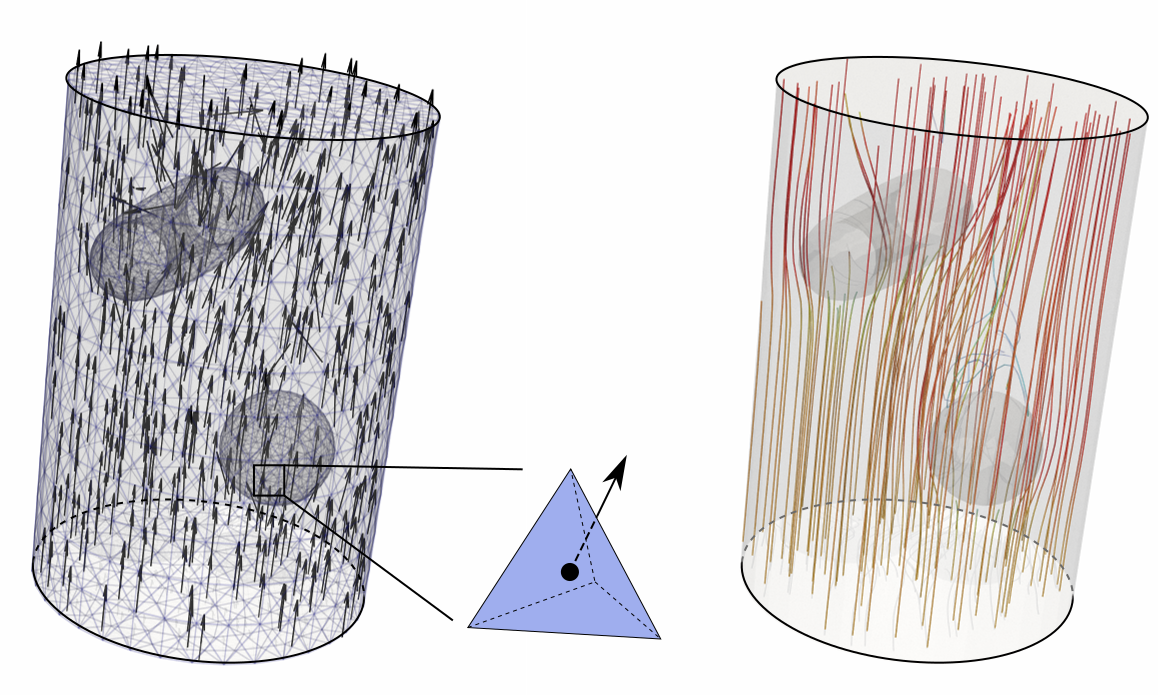}
	\caption{\label{fig:pvf}A piecewise constant vector field on a \texttt{Cylinder} tetrahedral mesh (left) and its streamline visualization using Runge-Kutta integration (right).}
\end{figure}

Let  $\surf$ be a simplicial solid, i.e a tetrahedral mesh in $\mathbb{R}^3$ with $n_v$ vertices, $n_e$ edges, $n_f$ faces, and $n_t$ tetrahedra. 
The set of simplices of dimension $d$ of $\surf$ is denoted by $\surf^{(d)}$. 
The number of interior simplices  and boundary simplices is denoted by $n_{ix}$ and $n_{\partial x}$, respectively, for $x\in\{v,e,f\}$.   A \textit{piecewise constant vector field} (PCVF) defined on $\surf$ is an element $X\in L^2(\surf,\mathbb{R}^3)$ such that $X_{|t}$ is constant in each tetrahedron $t$ (see figure~\ref{fig:pvf} for an example). We denote by $\vfsurf = \mathcal{X}(\mathcal{M}_h)$ the space of all PCVFs on $\surf$. 
The \textit{$L^2$-product} on $\vfsurf$ with respect to the Euclidean norm in $\mathbb{R}^3$ can then be written as the sum of weighted Euclidean scalar products 
\[
\left<X,Y\right>_{L^2} = \sum_{t\in \surf^{(3)}}\left<X_t,Y_t\right>_{\mathbb{R}^3}\cdot \mbox{vol}(t)
\]
where $\surf^{(3)}$ is the set of all tetrahedra of $\surf$ and $\mbox{vol}(t)$ is the Euclidean volume of a tetrahedron $t$. 

\subsection{Lagrange and Crouzeix-Raviart Spaces} 
We consider the following  function spaces define on $\surf$
\begin{align*}
\lagrange &= \lbrace\varphi:\surf \rightarrow \mathbb{R}^3\,\,|\,\, \varphi_{|t} \mbox{ is linear in } t \\
&\hspace{3.8cm} \mbox{ and } \varphi \mbox{ is globally continuous}\rbrace\\
\crouzeix &= \lbrace\psi:\surf \rightarrow \mathbb{R}^3\,\,|\,\, \psi_{|t} \mbox{ is linear in } t 
\\
&\hspace{2.2cm} \mbox{ and } \psi \mbox{ is continuous at facet 
	barycenters}\rbrace 
\end{align*}

The space $\lagrange$ is the space of linear \textit{Lagrange} elements on $\surf$. 
A basis of $\lagrange$ is given by the standard hat functions $\varphi_j$ satisfying $\varphi_j(v_i)=\delta_{ij}$  (Kronecker delta) for any  vertex $v_i\in\surf^{(0)}$. 
The space $\crouzeix$ is the space of \textit{Crouzeix-Raviart} elements. 
A basis of $\crouzeix$ is given  by a facet basis function $\psi_j$ with $\psi_j(b_i)=\delta_{ij}$, 
where $b_i$ is the barycenter of a (triangle) face $f_i \in \surf^{(2)}$. 
For the familiar case of simplicial surfaces, a  Crouzeix-Raviart element is linear over each triangle and continuous at edge midpoints whereas for tetrahedral meshes, they are linear over each tetrahedron and continuous at barycenters of the triangular  faces.

The boundary-constrained spaces $\lagrangebd \subset \lagrange$ and $\crouzeixbd \subset \crouzeix$ are defined by
\begin{eqnarray*}
	\lagrangebd &=&\lbrace \varphi \in \lagrange \,\,|\,\, \varphi(v_b)=0\mbox{ for all }v_b\in\partial \surf^{(0)} \rbrace\\
	\crouzeixbd &=&\lbrace \psi \in \crouzeix \,\,|\,\, \psi(b)=0\mbox{ for all barycenters }b\in\partial \surf^{(2)} \rbrace
\end{eqnarray*}
Note that elements in $\lagrangebd$ vanish in each boundary triangle whereas elements in $\crouzeixbd$ in general only vanish at the barycenters of boundary triangles. Our discrete decomposition in section \ref{sec:discrete_hodge_decomposition} does not involve the Lagrange function space. 
We only include it here to show the similarity between the two spaces. 
As pointed out in~\cite{poelke2017} the Langrage function space alone is not sufficient to obtain the correct dimensions in the discrete Hodge decomposition.

\subsection{Gradient Spaces} 
Given a piecewise linear function $\phi$ on $\surf$, its gradient $(\nabla \phi)_{|t}$ over a tetrahedron $t$ is defined as the gradient of $\phi_{|t}$, which can be expressed as a linear combination of the gradients of the basis functions restricted to $t$. 
The space of gradient vector fields is then defined as
\begin{eqnarray*}
	\nabla\lagrange &=&\lbrace \nabla\varphi \,\,|\,\, \varphi \in \lagrange \rbrace\\
	\nabla\crouzeix &=&\lbrace \nabla\psi  \,\,|\,\, \psi\in\crouzeix \rbrace
\end{eqnarray*}
By linearity, these gradients can be expressed as linear combinations of basis function gradients. Note that an element of the gradient space $\nabla\lagrangebd$ is then a vector field normal to $\partial\mathcal{M}$. 

\begin{figure}[!h]
	\centering
	\begin{tabular}{cc}
		\includegraphics[width=.4\linewidth]{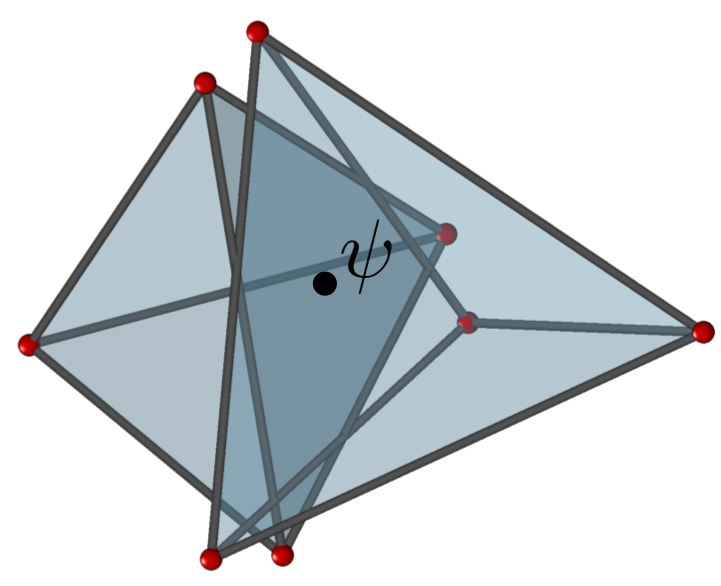}&\includegraphics[width=.3\linewidth]{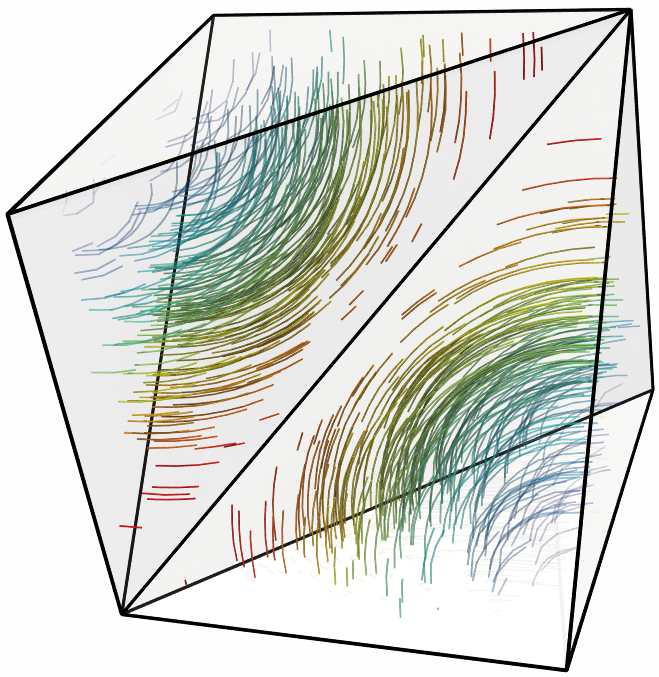}\\
		{\footnotesize (a) Crouzeix-Raviart function $\psi$} & {\footnotesize (b)  \textNedelec basis on an edge}
	\end{tabular}
	\caption{\label{fig:basis} A Crouzeix-Raviart basis function at a common face of two tetrahedra (left) and a \textNedelec basis function at a common edge of two tetrahedra (right).} 	
\end{figure}

\subsection{The \textNedelec Space}

The \textNedelec edge space $\nedelec$ is an edge-based finite element space used in the discretization of the Sobolev space  $\mbox{H(curl)}=\{v\in L^2(\smoothsurf)\,\,|\,\, \mbox{curl}(v)\in L^2(\smoothsurf)\}$. 
A convenient definition of this \textNedelec space makes use of differential forms and is based on its equivalence to the space of Whitney 1-forms~\cite{bossavit1988whitney}. 
In the language of vector proxies (cf.~\cite{nedelec_1980}), however, it can be defined in terms of the nodal basis functions $\varphi_i$ of $\mathcal{L}$,  interpreted as generalized barycentric coordinates on $\surf$.  
A basis of $\nedelec$ is then given by the set  $\{\varphi_i\nabla\varphi_j-\varphi_j\nabla\varphi_i\,\,|\,\,{e}_{ij}\in \surf^{(1)}\}$, indexed by the edges of $\surf$ (see figure~\ref{fig:basis} (b) for a visualization). 
Note that this explicit  definition requires a fixed orientation on the edges of $\surf$. 
A global \textNedelec basis function is non-zero over all tetrahedra sharing a given edge and zero everywhere else. 
In terms of vector fields, en element in $\nedelec$ is then linear in each tetrahedron, and tangentially continuous across interelement boundaries, i.e. the projection and restriction to the triangle from both sides between any two tetrahedra agrees. 
In analogy to the Lagrange and Crouzeix-Raviart spaces we denote by $\nedelec_0\subset \nedelec$ the subspace spanned by all \textNedelec basis functions associated to inner edges, and by 
$\curl(\nedelec_{(0)}) = \{\curl(\eta): \eta \in \nedelec_{(0)} \}$ the space of all curls of \textNedelec elements.

\section{Consistent Hodge-type decomposition}

We now present a consistent Hodge-type decomposition based on the previous ansatz spaces which is in one-to-one correspondence with the smooth setting. 
Using a straight-forward iterated projection procedure it is then easy to compute its decomposed components.
For details, we refer the reader to~\cite{poelke2017}.

\subsection{Orthogonal Relation and Harmonic Spaces}

The first fundamental relation between $\curl(\nedelec)$ and $\crouzeix$  is the insight that
\begin{align}
\curl(\nedelec) \perp \nabla\crouzeixbd\\
\curl(\nedelecbd) \perp \nabla\crouzeix. 
\end{align}
A proof of this relation can be found in (~\cite{poelke2017}, Lemma 3.3.2). 
Additionally, since $\nabla\crouzeixbd \subset \nabla\crouzeix$ and $\curl(\nedelecbd) \subset \curl(\nedelec)$, we have $\curl(\nedelecbd) \perp \nabla\crouzeixbd$. 
A visualization of elements in these spaces is provided in figure~\ref{fig:three}. 
Elements of $\curl(\nedelecbd)$ are vector fields which are vortices tangential to $\partial\smoothsurf_h$. 
Notice the importance of boundary conditions to assure orthogonality.

The $L^2$-orthogonal complements of the orthogonal sums $\nabla\crouzeixbd \oplus \curl(\nedelecbd)$, $\curl(\nedelec) \oplus \nabla\crouzeixbd$, and $\curl(\nedelecbd) \oplus \nabla\crouzeix$ within $\vfsurf$ are defined as the spaces of \textit{discrete harmonic fields} $\mathcal{H}_h$, \textit{discrete Neumann fields} $\mathcal{H}_{h,N}$ and \textit{discrete Dirichlet fields} $\mathcal{H}_{h,D}$, respectively. 
With these definitions, the spaces HG and HK correspond to the spaces $\neumann$  and  $\dirichlet$, respectively. 
They contain all curl- and divergence-free vector fields which depend solely on the topology of the shape. 
We have discrete de Rham isomorphisms $\neumann\cong H^2(\smoothsurf)$ and $\dirichlet\cong H^2(\smoothsurf,\partial \smoothsurf)$ between these spaces and the second cohomology and relative cohomology space, respectively, and we denote their dimensions by $h^2$ and $h_r^2$. 
In particular, the space $\dirichlet$ is nontrivial on shapes containing handles and tunnels whereas the space $\neumann$ is nontrivial for shapes with cavities.

\subsection{Discrete Hodge-Type Decompositions}
\label{sec:discrete_hodge_decomposition}
Now that we have defined all function spaces, we give several Hodge-type  decompositions involving these spaces. 
The decomposition starts with three fundamental terms. The additional subspaces are then derived by successive orthogonal projection.

\begin{lemma}[3-term fundamental decomposition]
	The space $\vfsurf$ admits the following $L^2$-orthogonal decompositions
	\begin{align}
	\vfsurf &= \curl(\nedelec)\oplus\nabla\crouzeixbd\oplus\neumann\label{eq:neumann}\\
	&= \curl(\nedelecbd)\oplus\nabla\crouzeix\oplus\dirichlet\label{eq:dirichlet}
	\end{align}
	with dim$\neumann=h^2$ and dim$\dirichlet=h^2_r$
\end{lemma}

\begin{figure}[!h]
	\centering
	\includegraphics[width=\linewidth]{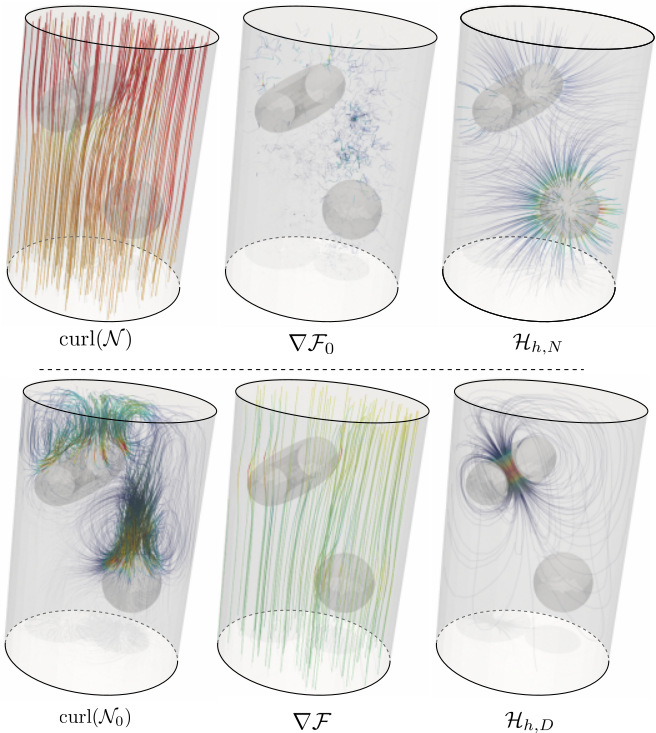}
	\caption{\label{fig:three} Visualization of each component of a 3-term Hodge decomposition of the vector field on the \texttt{Cylinder}.
		The first row is the fundamental Neumann, the second row the fundamental Dirichlet decomposition.}
\end{figure}

The proof of this Lemma follows directly from the previous observation. We call the first decomposition a \textit{fundamental Neumann decomposition} (FN) and the second decomposition a \textit{fundamental Dirichlet decomposition} (FD). If $\smoothsurf_h$ has a connected boundary, the spaces $\neumann$ and $\dirichlet$ vanish,  giving only two-component decompositions. 
Examples of vector fields representing each component of the decomposition of the \texttt{Cylinder} are depicted in figure~\ref{fig:three}.


\begin{theorem}[4-term, Hodge-Morrey-Friedrichs decompositon]
	The space $\vfsurf$ admits the following $L^2$-orthogonal decomposition
	\begin{align}
	\label{eq:fourterm1}
	\vfsurf &= \curl(\nedelecbd)\oplus\nabla\crouzeixbd\oplus\harmonic\cap\curl(\nedelec)\oplus\neumann\\\label{eq:fourterm2}
	&=\curl(\nedelecbd)\oplus\nabla\crouzeixbd\oplus\harmonic\cap\nabla\crouzeix\oplus\dirichlet.
	\end{align}	
\end{theorem}
\begin{figure}[!ht]
	\centering
	\includegraphics[width=0.7\linewidth]{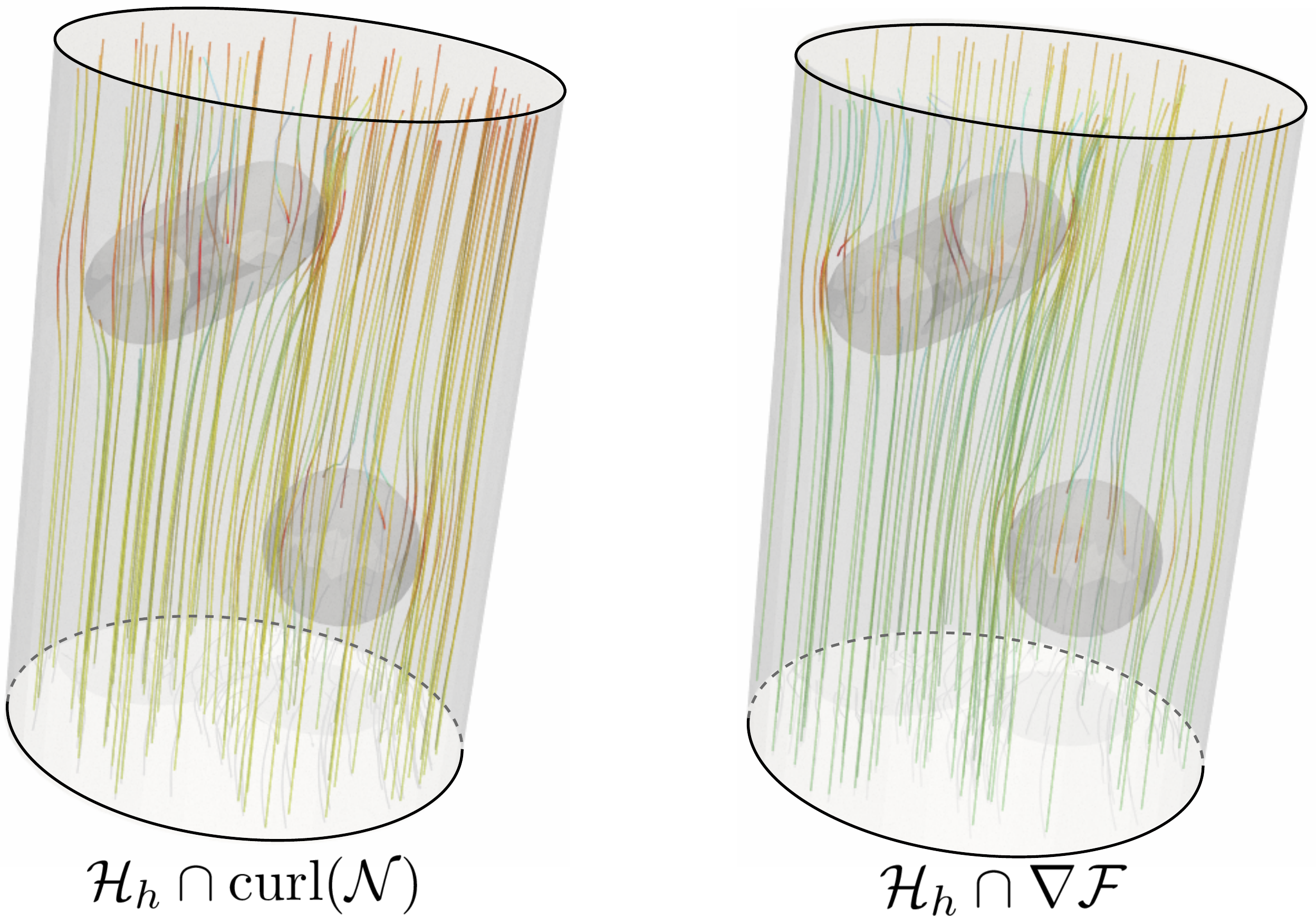}
	\caption{\label{fig:four} Hodge-Morrey-Friedrichs extra terms added to the fundamental decomposition.}
\end{figure}

The decomposition is a further extension of FN by expanding $\curl(\nedelec)$ into the orthogonal sum $\curl(\nedelec_0)\oplus\harmonic\cap\curl(\nedelec)$ and FD by expanding $\nabla\crouzeix$ into $\nabla\crouzeixbd\oplus\harmonic\cap\nabla\crouzeix$. 
These two spaces are very similar (see figure~\ref{fig:four} for the resulting vector field), the only difference is that $\harmonic\cap\curl(\nedelec)$ contains $\dirichlet$ while $\harmonic\cap\nabla\crouzeix$ contains $\neumann$. 
Finally we obtain a complete decomposition involving both harmonic spaces $\neumann$ and $\dirichlet$.

\begin{theorem}[5-term complete decomposition]
	The space $\vfsurf$ admits the following $L^2$-orthogonal decomposition
	\begin{equation}
	\vfsurf = \curl(\nedelecbd)\oplus\nabla\crouzeixbd\oplus\curl(\nedelec)\cap\nabla\crouzeix\oplus\neumann\oplus\dirichlet
	\end{equation}
\end{theorem}
The space $\curl(\nedelec)\cap\nabla\crouzeix$ is called the \textit{central harmonic component} obtained from equation~\eqref{eq:fourterm1} or ~\eqref{eq:fourterm2} by further decomposing the space $\harmonic$.  
It has dimension $n_{\partial f}-h^2-1$ which clearly depends on the mesh resolution and is not topologically invariant. 	In particular, it increases under uniform refinement of the mesh.  
This is in accordance with the smooth setting, as it corresponds to the infinite-dimensional space CG of curly gradients. 


\subsection{Iterative  $L^2$-Projection}
We use a standard iterative $L^2$-projection to compute a decomposition. 
Starting with a 3-term decomposition and depending on the topology of the simplicial manifold, the 4-term and  5-term decompositions are derived.

The basic idea of the iterative $L^2$-projection is to subtract a component from the previous vector field after projection. Suppose we have a vector field $X=X_0\in\vfsurf$ that we want to decompose. Then  $X_{i+1} = X_i - X_{\mathcal{Y}}$, where $X_{\mathcal{Y}} = \texttt{Project}(X_i, \mathcal{Y})$ and $\mathcal{Y}$ is one of the spaces in the discrete decomposition of interest. 
In the following we illustrate this procedure for the FD decomposition. 
The FN decomposition can be computed analogously.

Starting with the input vector field $X$, it is projected first onto $\mbox{curl}(\nedelecbd)$. The difference $X-\xnedelecbd$ is then projected onto $\nabla\crouzeix$. 
Finally, the remainder defines the harmonic Dirichlet field $\xdirichlet$ which is naturally a tangential vector field.

To get the four term decomposition it is enough to project $\xcourzeix$ onto $\nabla\crouzeixbd$ which means extracting the component of the gradient normal to the boundary. The remainder is then the new component $\xharmcoex$. The five term decomposition is finally obtained by projecting $\xharmcoex$ onto $\curl(\nedelec)$. The Neumann harmonic field is derived as $\xneumann=\xharmcoex - \xcentre$.

Note that throughout the different decompositions, all  new components are subspaces of the space $\nabla\crouzeix$ or $\curl(\nedelec)$. 
In fact we even have $\vfsurf = \curl(\nedelec) + \nabla\crouzeix$, although this sum is neither orthogonal nor direct, cf. \cite[Cor.~3.3.10]{poelke2017}. 
If the $L^2$-norms of these components are initially small, all subsequent components are small, too, or even negligible depending on the nature of the vector fields. 
However, in some applications the geometry of the streamline paths is more important than the local magnitude of the vectors.

\subsection{Computing the Projections} The $L^2$-projection \verb|Project| follows the traditional Galerkin method. Let $\mathcal{U}$ and $\mathcal{V}$ be two spaces of vector fields with the associated basis $\mathcal{B}_{\mathcal{U}} = \{\mu_i\}_{i=1,\ldots,n}$ and $\mathcal{B}_{\mathcal{V}} = \{\nu_i\}_{i=1,\ldots,m}$. Suppose $\hat{u}$ is the $L^2$-projection of $u$ in $\mathcal{V}$, then 
\begin{equation*}
\hat{u} = \sum_{i=1}^{m}\hat{u}_i\nu_i\mbox{ and }\left<u-\hat{u},\nu_j\right>=0,\,\,\forall j=1,\ldots,m.
\end{equation*}
Writing $u=\sum_i u_i \mu_i$ and rearranging the terms yields the following linear equality

\begin{equation*}
\sum_{i=1}^{m}\hat{u}_i\left<\nu_i,\nu_j\right> 
= \sum_{i=1}^n u_i\left<\mu_i,\nu_j\right>,\,\,\forall j=1,\ldots,m.
\end{equation*}
In matrix form, it becomes 
$$M_{\scriptscriptstyle\mathcal{V}}\hat{\mathbf{u}} = P_{\scriptscriptstyle \mathcal{U},\mathcal{V}}\mathbf{u},$$ 
where $M_{\scriptscriptstyle\mathcal{V}}$ is the mass matrix of $\mathcal{V}$ and $P_{\mathcal{U},\mathcal{V}}$ is the projection matrix from the space $\mathcal{U}$ to $\mathcal{V}$. 


\section{Numerical Evaluation}
In the following we evaluate the stability and scalability of our approach using blood flow vector fields acquired from a standard procedure in medical computational fluid dynamic. 

\subsection{Implementation Details}

We use the \fenics library \cite{fenicsweb} in version 2018.1.0,  combined with Python to compute the discrete Hodge decomposition. 
\fenics provides a domain specific language interface as well as many classical function spaces for finite element problem so that the 
actual implementation closely resembles pseudo-code examples, cf. figures~\ref{fig:pseudo1} and~\ref{fig:pseudo2}.  
We state the projection routine used  in the case where the differential operator is the curl with a Dirichlet boundary condition. The other operator and boundary condition can be handled similarly. 

The projection of a  vector field $X$ onto the space $\curl(\nedelecbd)$ in Python can be done following the code snippet in figure~\ref{fig:pseudo1}.

\begin{figure}[!h]
	\begin{adjustbox}{width=.9\columnwidth,center}
		\begin{tabular}{|l|}
			\hline
			{\footnotesize{\color{brown} \verb|# Mesh and vector field in dolfin format|}}\\
			\verb|M_h = Mesh(|{\color{red}\verb|"path_to_mesh"|}\verb|)|\\
			\verb|vf =| \verb|VectorFunctionSpace(M_h,|{\color{red}\verb|`DG'|}\verb|,0)|\\
			\verb|X_h = Function(vf,|{\color{red}\verb|"path_to_vf"|}\verb|)|\\\\
			{\footnotesize{\color{brown}\verb|# Nedelec space|}}\\
			\verb|Nedelec = FunctionSpace(M_h,|{\color{red}\verb|`N1curl'|}\verb|,1)|\\\\
			{\footnotesize{\color{brown}\verb|# Define variational problem|}}\\
			\verb|u = TrialFunction(Nedelec)|\\
			\verb|v = TestFunction(Nedelec)|\\
			\verb|M = inner(curl(u),curl(v))*dx|\\
			\verb|P = inner(X_h,curl(v))*dx|\\
			\verb|c = Constant((0.,0.,0.))|\\
			\verb|bc = DirichletBC(Nedelec,c, DomainBoundary())|\\\\
			{\footnotesize{\color{brown}\verb|# Compute solution|}}\\
			\verb|v = Function(V)|\\
			\verb|solve(M == P, v, bc)|\\\\
			{\footnotesize{\color{brown}\verb|# Change basis|}}\\
			\verb|pr_vf = project(curl(v),X_h)|\\\\
			\hline    
		\end{tabular}
	\end{adjustbox}
	\caption{Projecting a vector field $X$ onto space curl($\nedelec_0$) in \fenics}
	\label{fig:pseudo1}
\end{figure}

In this code, the input mesh is represented by \verb|M_h| and the vector field by \verb|X_h|. 
The function space \verb|`DG'| stands for discontinuous Galerkin elements, here of order $0$, which is just a standard Euclidean representation of a single constant vector per tetrahedron. 
The space \verb|`N1curl'| is the lowest order \textNedelec edge-element space of the first kind. 
The rest of the code is the standard variational formulation including the curl operator. 
The quantity \verb|v| represents a vector potential whereas \verb|pr_vf| is the actual solution of the projection, represented as a PCVF.

To compute the projection onto $\nabla\crouzeix$, the differential operator \verb|curl| is  replaced by \verb|grad|, 
the space \verb|N1curl| is replaced by \verb|CR|, which stands for Crouzeix-Raviart elements, and the boundary condition \texttt{bc} is removed. 
All functions given in the above code snippet are available in \fenics. 
What remains is the computation of the difference of two PCVFs, which is given in figure~\ref{fig:pseudo2}.

\begin{figure}[!h]
	\begin{adjustbox}{width=.8\columnwidth,center}
		\begin{tabular}{|l|}
			\hline
			{\footnotesize{\color{brown} \verb|# Input vector field X_h and Y_h|}}\\
			{\footnotesize{\color{brown} \verb|# Function space which should be the same|}}\\
			\verb|V = X_h.function_space()|\\\\
			{\footnotesize{\color{brown} \verb|# Define difference vector|}}\\
			\verb|sub = Function(V)|\\\\
			{\footnotesize{\color{brown} \verb|# Assign entry|}}\\
			{\footnotesize\verb|sub.vector()[:] = x.vector().get_local()|}\\
			\verb|              - y.vector().get_local()|\\\\
			\hline    
		\end{tabular}
	\end{adjustbox}
	\caption{\label{fig:pseudo2}Subtracting two vector fields of the same function space in \fenics}
\end{figure}


Because of the $L^2$-orthogonality of the decomposition, the Pythagorian formula gives a full quantification for each component: 
the squared $L^2$-norm of the input vector field is equal to the some of the squared $L^2$-norms of each components. 

We use \paraview \cite{paraviewweb} for the streamline visualization. 
The individual streamlines are color-coded with opacity control according to the local magnitude of the vectors in each cell. 
In our color scheme white (resp. transparent) encodes low magnitude while dark red (resp. opaque) encodes high magnitudes.

\subsection{Validating Orthogonality}

To check numerically the $L^2$-orthogonality of the discrete decomposition, we run the algorithm on the example fields given in figure~\ref{fig:example}. 
We expect that each vector field $X_0$, \dots , $X_4$ shows up in its respective space in the discrete setting. 
Moreover, we compute the decomposition for the superposition $X_{012} = X_0+X_1+X_2$. 
If the $L^2$-norm of a resulting component is below $10^{-10}$  we  consider it to be zero. 
The results are shown in table~\ref{tab:decomp}. 
Note how the values identify the dominating components in the decomposition, both for the ``pure'' fields $X_0$, \dots, $X_4$ as well as for the vector field sum $X_{012}$.
Of course there are approximation errors in the discrete setting, stemming from the interpolation of the smooth fields, the approximation of the 
smooth geometry and the numerical errors in the solver of the linear system. 
As we will see in the next section, the mesh resolution has a strong impact on the precision. 

\begin{table}[!h]
	\centering
	\begin{tabular}{|c|l|l|l|l|l|l|}
		\hline
		&&&&&&\\
		& $\vfsurf$ & $\curl(\nedelec_0)$ & $\nabla\crouzeixbd$ & $\centre$ & $\neumann$ & $\dirichlet$\\ \hline
		&&&&&&\\
		$X_0$             & 1.603  & 1.602  & 8-e04  & 3-e04  & 0.  & 0. \\
		$X_1$     & 2.403  & 0.0    & 2.403    & 2-e04  & 0. & 0.\\
		$X_2$   & 3.065  & 0.  & 0.    & 3.065  & 0.     & 0.\\
		$X_3$             & 37.68  & 9-e04  & 1.01  & 0.006  & 36.7     & 0.\\
		$X_4$            & 10.04    & 1.7-e05    & 0.028    & 0.002  & 0.       & 10.\\
		$X_{012}$          & 7.071    & 1.602    & 2.404    & 3.064    & 0.     & 0.\\
		&&&&&&\\
		\hline
	\end{tabular}
	\caption{The $L^2$-norms square of each component of the Hodge decomposition of the vector fields defined in equation~\eqref{eq:vfsmooth}.}
	\label{tab:decomp}
\end{table}

\begin{landscape}
	\begin{figure*}
	\includegraphics[width=\linewidth]{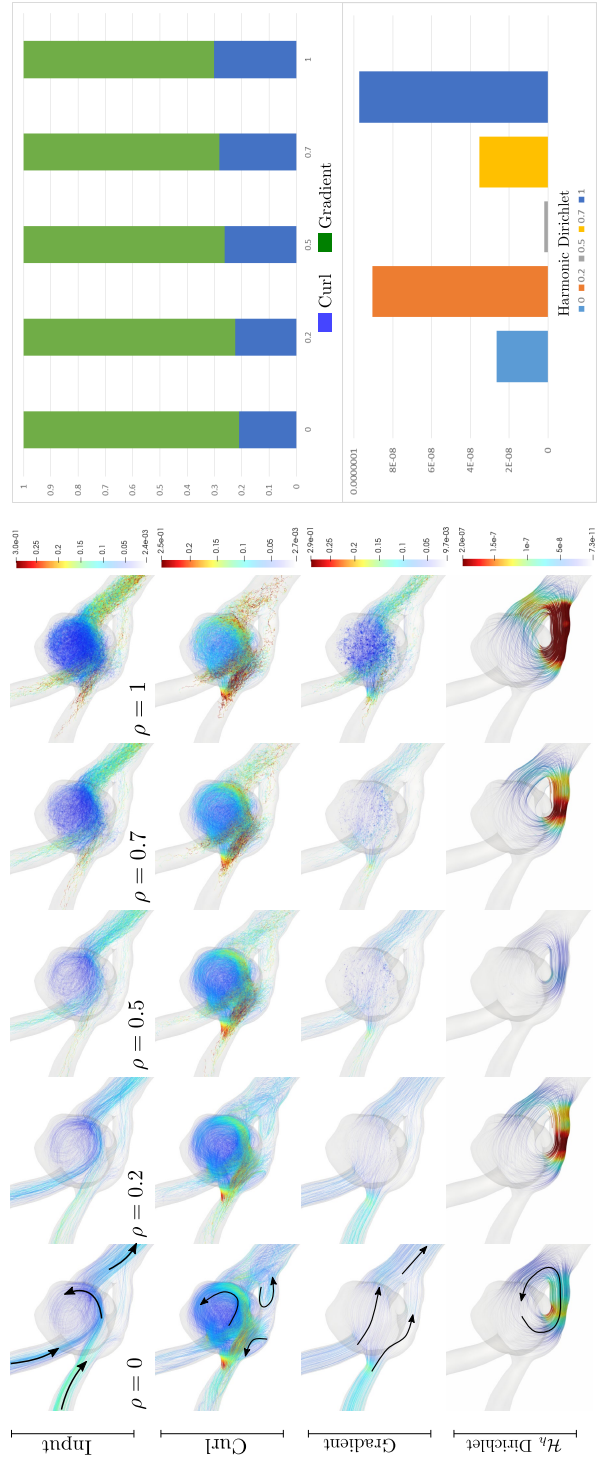}
	\caption{\label{fig:noise} FD decomposition of a progressive Gaussian noised vector field with factor $\rho = 0, .5, .7, 1.$ on a fenestrated aneurysm. The strong features of each component are still visually identifiable even for large noise frequency. The geometry of the harmonic Dirichlet streamline is invariant under noise. The diagram shows the normalized $L^2$-norms of each component.}
	\end{figure*}
\end{landscape}

\subsection{Generating 3D Vector Fields}
\label{subsec:generating_3d_vfs}
To generate 3D vector fields,  we use simulated blood flow following  the same methodology as described in~\cite{razafindrazaka2019}. 
In summary, the geometries are first acquired by a routinely used imaging modality including segmentation and smoothing (using ZIBAmira 2015.28, Zuse Institute Berlin,  Germany). The blood is then modeled as a non-Newtonian fluid with a constant density of $1050 \mathrm{kgm}^{-3}$ and a Carreau-Yasuda viscosity model~\cite{karimi2014}. 
All simulations were performed as steady-state simulations of the peak-systolic flow using STAR-CCM+ (v. 12.06, Siemens PLM Software, Plano, USA). The result is then a tetrahedral mesh with a piecewise linear vector field. 
We use the FD decomposition whose curl component is tangential to the boundary surface, as the components of this decomposition better capture the features of the flow, such as vortices and laminar flows, in contrast to the FN decomposition (see figure~\ref{fig:three}). 

\subsection{Stability Under Varying Mesh Resolution}
\begin{figure}[!h]
	\includegraphics[width=\linewidth]{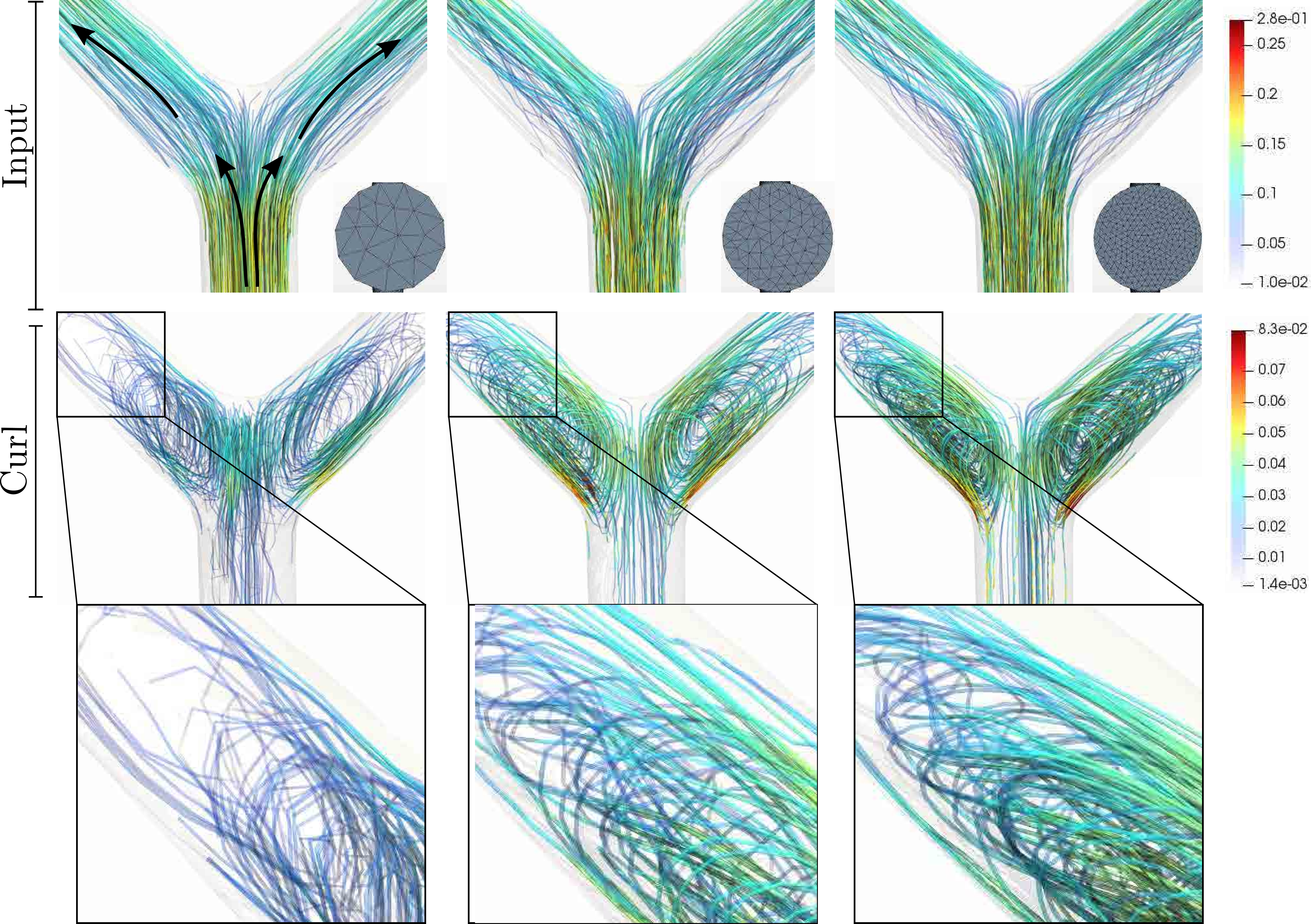}
	\caption{\label{fig:res} Curl component analysis of a vector field under mesh refinement: 4k, 61k, and 161k. 
		The effect of the piecewise linear approximation is visible in the geometry of the pattern.}
\end{figure}

In this section, we analyze the stability of the decomposition under mesh refinement. 
The most problematic regions of the mesh are usually those close to the boundary surface where the tetrahedra need to be very fine in order to approximate the curvature of the shape. 
In our experiments we use three resolutions (4k, 61k, and 161k tetrahedra) of a model and perform the decomposition for  simulated blood flow on these models. 
Figure~\ref{fig:res} gives the results on a simple model. 
It shows only the curl component since the gradient component is very similar to the input field on this model.   
The squared $L^2$-norms scaled by $10^4$ are given in table~\ref{tab:decomp_mesh}. 
There is a large difference between the norms on the low resolution models, which is no surprise since the approximation quality is lower. Note, however, that all numerical studies assume a preliminary mesh independence study for the parameter of the interest. In cases of biological flow problems, this is usually a wall shear stress, which requires a very high resolution mesh.
\begin{table}[!h]
	\centering
	\begin{tabular}{|c|l|l|l|}
		\hline
		&&&\\
		& 4k & 61k & 161k \\ \hline
		&&&\\
		$X_h$             & 0.547  & 0.509  & 0.509 \\
		$\curl(\nedelec_0)$     & 0.005  & 0.014    & 0.019 \\
		$\nabla\crouzeix$   & 0.541  & 0.494  & 0.490\\    
		&&&\\
		\hline
	\end{tabular}
	\caption{The squared $L^2$-norms  of each component of the Hodge decomposition of the vector fields for different mesh resolutions 4k, 61k, and 161k.}
	\label{tab:decomp_mesh}
\end{table}

\subsection{Noisy Vector Field}
To evaluate the robustness of the method, a simulated vector field using the procedure in section~\ref{subsec:generating_3d_vfs} on a fenestrated aneurysm is used. 
Again we use the FD decomposition to identify the important change in the flow. 
The genus of the model causes the space of Dirichlet fields to be nontrivial this time. 
Gaussian noises are added increasingly to the input vector fields following the method proposed in \cite{random1998}. 
More precisely, for each vector  $\mathbf{v}$ of a vector field, we add the value  $\rho\|\mathbf{v}\|\sigma_{\mathbf{v}}$  to $\mathbf{v}$. Here $\rho$ is a noise factor. The values  0.2, 0.5, 0.7, and 1 are chosen for the experiments. 
For each of these noised vector fields, we apply the decomposition and compare the velocity magnitude with respect to the original noise-free vector field.

The corresponding noisy vector fields and the decomposition are shown in figure~\ref{fig:noise}. The first row represents the progressive noised vector fields from $\rho=0$ to $1$. The second to the last row are the components of the decomposition. Notice that the Dirichlet harmonic flow has a very low contribution to the field but remains geometrically invariant under the noise. It has a maximal $L^2$-square norm of 1e-07. For the other component, some of the strong features (with high magnitude) of each component are still recovered in the decomposition. The low magnitude features are converted into noise. Notice the similarity of the  vortices' pattern in the tumor for both the original fields and the noisy fields. 
Comparing the normalized norms of each component, we notice a slight increase in the curl component  relative to the noise magnitude. 

The stability of the Hodge decomposition analysis regarding noise could be very helpful for an analyis of all experimentally measured vector fields used as the visualzation technique solely (e.g. MRI in the medicine) or as a validation approach for CFD.

Note that our method cannot handle measured velocity fields acquired from other procedures such as MRI without further pre-processing. 
One could eventually map the voxel-based MRI information onto a finite element grid but this procedure is delicate and is out of the scope of this paper.  
However, the noise analysis shows that the decomposition could still identify important features of the field. 

\subsection{Anatomical Change}

We use the Hodge decomposition to identify the blood flow change in a patient stenosis before and after intervention. Coarctation of the aorta causes abnormal blood flow conditions in the aorta behind the narrowing associated with a formation of a turbulent jet and a recirculation zone resulting in a high pressure drop through the narrowing. The medical intervention is then intended to re-open the narrow segment by using surgical or catheter-based procedures to allow for a normal blood flow. Before the intervention, the stenosis is scanned through CT and blood flow is simulated using a constant input velocity profile. The same process is repeated after the intervention. To obtain a continuous deformation of the pre to the post operative model, we use a hyperelastic interpolation as introduced in \cite{Hildebrandt2011b}. The idea is to mimick the formation of a stenosis and the corresponding deterioration or improvement in blood flow. 

\begin{figure}[!h]
	\includegraphics[width=\linewidth]{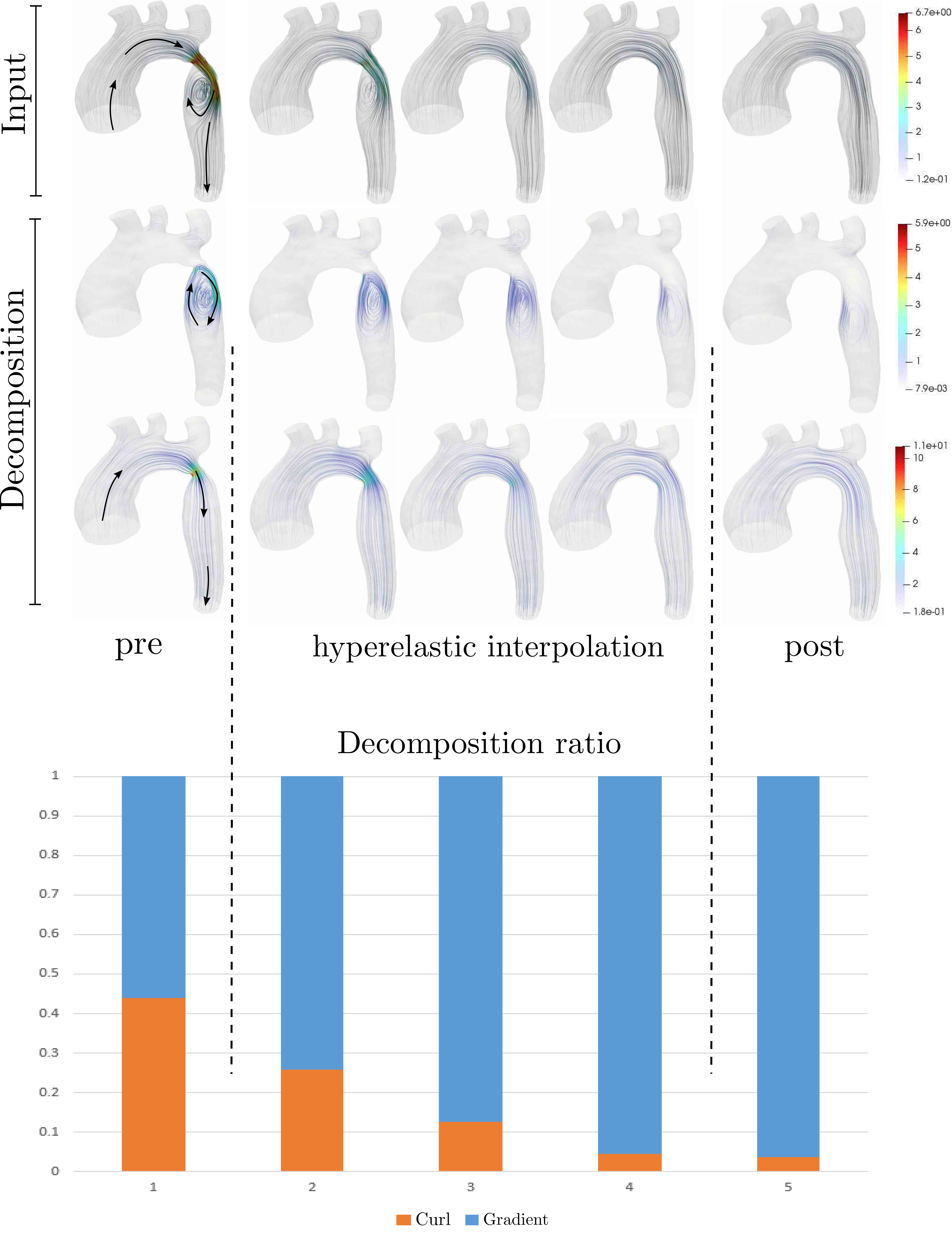}
	\caption{\label{fig:anatom} Decomposition of a simulated blood flow of a patient before and after intervention. A hyperelastic interpolation following \cite{Hildebrandt2011b} is applied to see the evolution of each component of the decomposition.}
\end{figure}

The results are given in figure \ref{fig:anatom} for three in-between frame models. 
The first observation is the reduction in the curl component after intervention deduced from the component ratio analysis. 
Geometrically, the large vortices are mainly located behind the stenosis and are away from the inlet boundary. 
Notice also the improvement in the gradient flow which becomes more uniform in terms of velocity magnitude. 
Using other boundary conditions such as in \cite{Zhao2019} may give a completely different ratio, but as the vortices are located away from the inlet, we believe that having a free boundary condition at the inlet would not affect the results very much.

\subsection{Boundary Conditions in CFD}
Boundary conditions are very important to obtain an accurate simulation of a patient specific physiological blood flow. 
Since only a small segment of the patient's vessel is analyzed, two methods are popularly used to set the boundary conditions at the inlet velocity profile: vector information based on magnetic resonance imaging (MRI) measurements or simplified profiles including constant (plug), parabolic or alternatively Womersley input velocity profile with a MRI-measured or literature-based inlet flow rate~\cite{gallo2012,ladisa2011,vergara2011}. 
On the basis of  wall shear stress on surfaces it was shown in \cite{razafindrazaka2019}  that the difference of both methods is reflected in the harmonic flows. 
\begin{figure}[!ht]
	\centering
	\includegraphics[width=\linewidth]{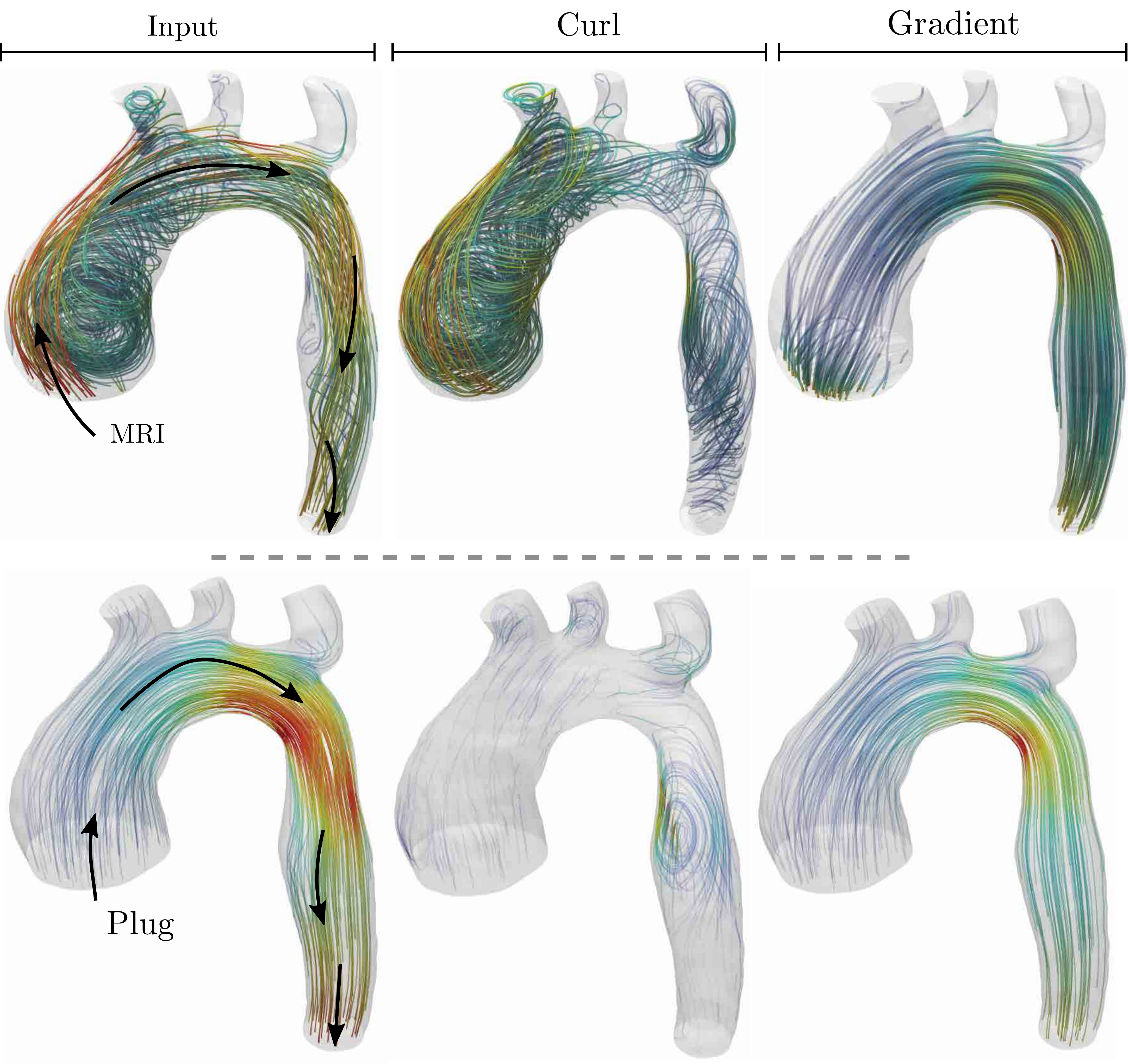}
	\caption{Comparison of FD decompositions  measured velocity profile from MRI scan and constant input as inlet boundary condition.}
	\label{fig:mri_vs_plug}
\end{figure}


Figure~\ref{fig:mri_vs_plug} shows the resulting decomposition using MRI-based inlet velocity profile after intervention on the same patient as in the anatomical analysis in the previous section. Notice the huge amount of curl component close to the inlet surface. These partitions can be seen through the decomposition where the gradient field is less dominant in this region but more present on other region. This is no surprise since the MRI inlet profile is quite complex inducing a more turbulent flow. Comparing the normalized $L^2$-norms between MRI and plug boundary condition, we have 41\% vs 9\% for the curl component and 59\% vs 91\% for the gradient component.

\section{Limitation and Future Work}
Our implementation is based on the formulation of~\cite{poelke2017} which interprets vector proxies as 2-forms. 
The $L^2$-product with gradients of Crouzeix-Raviart elements gives a natural measure for ``normal continuity'' of vector fields across interelement boundaries. Normal continuity is necessary  to define a weak divergence. 
Doing the same for 1-forms would be more attractive since it would avoid swapping the  boundary conditions, resulting in a direct discretization of the decompositions for vector fields representing 1-forms. 
It  would require a finite element ansatz space such that the $L^2$-orthogonal complement of the curl  of its elements is the space of ``tangentially continuous'' PCVFs. 
We are not aware of a convenient definition for such an ansatz space, though.

Our formulation and implementation of the 3D Hodge decomposition leaves room for improvement in various directions. 	
One of them is running time. 
Medical data coming from CFD simulations exhibits a huge number of tetrahedra required for high accuracy in numerical simulations. 
A multiresolution approach e.g. by subspace approximation could speed up the computation considerably. 
One could also simplify the tetrahedral mesh by maintaining all important properties of the field ~\cite{platis2004}.


\section{Conclusion}
In this article we reviewed a consistent 3D Hodge decomposition with detailed implementation and practical evaluation.	
Using a modern finite element library, it can be implemented in a few lines of code. 
Our discretization of the Hodge decomposition was for the first time applied to analyze the 3D velocity vector fields of simulated patient-specific blood flows as a pilot study.  
It is able to distinguish between pathological and physiologic blood flows, and to characterize the impact of inflow boundary conditions as well as the impact of a treatment procedure. 

\section*{Acknowledgment}
The authors would like to thank Pavel Yevtushenko for the geometries with simulated blood flow.  This project  is funded by the Deutsche Forschungsgemeinschaft (DFG, German Research Foundation) under Germany’s Excellence Strategy – The Berlin Mathematics Research Center MATH+ (EXC-2046/1, project ID: 390685689) and partly funded by the Bundesministerium für Bildung und Forschung (Project VIP+, DSSMitral).

\bibliographystyle{unsrt}
\bibliography{phd_references}

\end{document}